\documentstyle[amscd,amssymb,verbatim,epsf]{amsart}
\begin{document}

\theoremstyle{plain}
\newtheorem{Thm}{Theorem}
\newtheorem{Cor}{Corollary}
\newtheorem{Con}{Conjecture}
\newtheorem{Main}{Main Theorem}

\newtheorem{Lem}{Lemma}
\newtheorem{Prop}{Proposition}

\theoremstyle{definition}
\newtheorem{Def}{Definition}
\newtheorem{Note}{Note}

\theoremstyle{remark}
\newtheorem{notation}{Notation}
\renewcommand{\thenotation}{}

\errorcontextlines=0
\numberwithin{equation}{section}
\renewcommand{\rm}{\normalshape}%

\title[$p$-quasi-Cauchy sequences]
   {$p$-quasi-Cauchy sequences}
\author{H\"usey\.{I}n \c{C}akall\i \\Maltepe University,TR 34857, Maltepe, Istanbul, Turkey Phone:(+90216)6261050 ext:2248, \;  fax:(+90216)6261113}
\address{H\"Usey\.{I}n \c{C}akall\i\\
          Department of Mathematics, Maltepe University, Marmara E\u{g}\.{I}t\.{I}m K\"oy\"u, TR 34857, Maltepe, \.{I}stanbul-Turkey \; \; \; \; \; Phone:(+90216)6261050 ext:2248, fax:(+90216)6261113} \email{hcakalli@@maltepe.edu.tr; hcakalli@@gmail.com}
\date{\today}

\begin{abstract}

In this paper we generalize the concept of a quasi-Cauchy sequence to a concept of a $p$-quasi-Cauchy sequence for any fixed positive integer $p$. For $p=1$ we obtain some earlier existing results as a special case. We obtain some interesting theorems related to $p$-quasi-Cauchy continuity, $G$-sequential continuity, slowly oscillating continuity, and uniform continuity. It turns out that a function $f$ defined on an interval is uniformly continuous if and only if there exists a positive integer $p$ such that $f$ preserves $p$-quasi-Cauchy sequences where a sequence $(x_{n})$ is called $p$-quasi-Cauchy if $(x_{n+p}-x_{n})_{n=1}^{\infty}$ is a null sequence.
\end{abstract}

\maketitle

\section{Introduction}
\normalfont
The concepts of continuity, compactness, and connectedness and any other concept involving these concepts play a very important role not only in pure mathematics but also in other branches of science involving mathematics especially in computer science, combinatorics, information theory, biological science, geographic information systems, population modeling, and motion planning in robotics.

Recently, in \cite{CakalliForwardcontinuity}, a concept of quasi-Cauchy continuity,  and a concept of quasi-Cauchy compactness have been introduced in the senses that a real function is called quasi-Cauchy continuous if $\lim_{n\rightarrow\infty}\Delta f(x_{n})=0$ whenever $\lim_{n\rightarrow\infty}\Delta x_{n}=0$, and a subset $E$ of $\textbf{R}$ is called quasi-Cauchy compact if whenever $(x_{n})$ is a sequence of points in $E$ there is a subsequence $(y_{k})=(x_{n_{k}})$ of $(x_{n})$ with $\lim_{k\rightarrow \infty} \Delta y_{k}=0$ where $\Delta y_{k}=y_{k+1}-y_{k}$. We note that forward continuity and forward compactness terms were used in place of the terms quasi-Cauchy continuity and quasi-Cauchy compactness in \cite{CakalliForwardcontinuity}, respectively.
What happens if we replace $\Delta y_{n}=y_{n+1}-y_{n}$ by $\Delta y_{n}=y_{n+2}-y_{n}$. As a matter of fact we could replace $\Delta y_{n}=y_{n+1}-y_{n}$ by $\Delta y_{n}=y_{n+p}-y_{n}$ for any positive integer $p$.

The purpose of this paper is to generalize the concept of ward continuity to a concept of $p$-quasi-Cauchy continuity for any fixed positive integer $p$, and investigate relations among such kinds of continuities, slowly oscillating continuity, uniform continuity, $G$-sequential continuity, and ordinary continuity.


\maketitle

\section{\lowercase{Preliminaries}}
Now we give some definitions and notation which will be needed throughout the paper. $\textbf{N}$ will denote the set of all positive integers. We will use boldface letters $\bf{x}$, $\bf{y}$, $\bf{z}$, ... for sequences $\textbf{x}=(x_{n})$, $\textbf{y}=(y_{n})$, $\textbf{z}=(z_{n})$, ... of terms in $\textbf{R}$. $c$, and $\Delta$ will denote the set of all convergent sequences, and the set of all quasi-Cauchy sequences of points in $\textbf{R}$ where a sequence $\textbf{x}=(x_{n})$ is called quasi-Cauchy if $\lim_{n\rightarrow\infty}\Delta x_{n}=0$ which was given in \cite{CakalliForwardcontinuity} as being called forward convergent to $0$ (see also \cite{BurtonColeman}).
It is known that a sequence $(x_{n})$ of points in $\textbf{R}$, the set of real numbers, is slowly oscillating if $$
\lim_{\lambda \rightarrow 1^{+}}\overline{\lim}_{n}\max _{n+1\leq
k\leq [\lambda n]} |
  x_{k}  -x_{n} | =0
$$ where $[\lambda n]$ denotes the integer part of $\lambda n$. This is equivalent to the following if $x_{m}-x_{n}\rightarrow 0$ whenever $1\leq \frac{m}{n}\rightarrow 1$ as, $m,n\rightarrow \infty$.
Using $\varepsilon>0$ and $\delta$ this is also equivalent to the case when for any given $\varepsilon>0$, there exists $\delta=\delta (\varepsilon) >0$ and $N=N(\varepsilon)$ such that $|x_{m}-x_{n}|<\varepsilon$ if $n\geq N(\varepsilon)$ and $n\leq m \leq (1+\delta)n$ (see \cite{CakalliSlowlyoscillatingcontinuity}).
Any Cauchy sequence is slowly oscillating, and any slowly oscillating sequence is quasi-Cauchy. But the converses are not always true. For example, the sequences $(\sum^{\infty}_{k=1}\frac{1}{n})$, (ln n), (ln ln n), and combinations like that are slowly oscillating, but not Cauchy. The sequence $(\sum^{k=n}_{k=1}(\frac{1}{k})(\sum^{j=k}_{j=1}\frac{1}{j}))$ is quasi-Cauchy, but not slowly oscillating (see also \cite{FDikMDikandCanak}, \cite{DikandCanak}, \cite{Vallin} and \cite{VallinandCakalli}).

The idea of statistical convergence was formerly given under the name "almost convergence" by Zygmund in the first edition of his celebrated monograph published in Warsaw in 1935 \cite{ZygmundTrigonometricseries}. The concept was formally introduced by Fast \cite{Fast} and later was reintroduced by Schoenberg \cite{SchoenbergTheintegrabilityofcertainfunctions}, and also independently by Buck \cite{BuckGeneralizedasymptoticdensity}. Although statistical convergence was introduced over nearly the last seventy years, it has become an active area of research for twenty years. This concept has been applied in various areas such as number theory \cite{ErdosSurlesdensitesde}, measure theory \cite{MillerAmeasuretheoresubsequencecharacterizationofstatisticalconvergence}, \cite{ChengLinLanLiuMeasuretheoryofstatisticalconvergence}, trigonometric series \cite{ZygmundTrigonometricseries}, summability theory \cite{FreedmanandSemberDensitiesandsummability}, locally convex spaces \cite{MaddoxStatisticalconvergenceinlocallyconvex}, in the study of strong integral summability \cite{ConnorandSwardsonStrongintegralsummabilityandstonecompactification}, turnpike theory \cite{MakarovLevinRubinovMathematicalEconomicTheory}, \cite{MckenzieTurnpiketheory}, \cite{PehlivanandMamedovStatisticalclusterpointsandturnpike}, Banach spaces \cite{ConnorGanichevandKadetsAcharacterizationofBanachspaceswithseparabledualsviaweakstatisticalconvergence}, and metrizable topological groups \cite{Cakallilacunarystatisticalconvergenceintopgroups}, and topological spaces \cite{MaioKocinac}, \cite{CakalliandKhan}. It should be also mentioned that the notion of statistical convergence has been considered, in other contexts, by several people like R.A. Bernstein, Z. Frolik, etc. The concept of statistical convergence is a generalization of the usual notion of convergence that, for real-valued sequences, parallels the usual theory of convergence. A sequence $(x _{k})$ of points in $\textbf{R}$ is called statistically convergent to an element $\ell$ of $\textbf{R}$  if for each
$\varepsilon$
\[
\lim_{n\rightarrow\infty}\frac{1}{n}|\{k\leq n: |x_{k}-\ell|\geq{\varepsilon}\}|=0,
\] and this is denoted by $st-\lim_{k\rightarrow\infty}x _{k}=\ell$ \cite{Fridy} (see also \cite{CakalliAstudyonstatisticalconvergence}). $S(\textbf{R})$ will denote the set of statistically convergent sequences of points in $\textbf{R}$.

A real sequence $(x_{k})$ is called lacunary statistically convergent to an element $\ell$ of $\textbf{R}$ if
\[
\lim_{r\rightarrow\infty}\frac{1}{h_{r}} |\{k\in I_{r}: |x_{k}-\ell|\geq{\epsilon} \}|=0,
\]
for every $\epsilon>0$ where $I_{r}=(k_{r-1},k_{r}]$ and $k_{0}=0$,
$h_{r}:k_{r}-k_{r-1}\rightarrow \infty$ as $r\rightarrow\infty$ and $\theta=(k_{r})$ is an increasing sequence of positive integers and, and this is denoted by $S_{\theta}-\lim_{n\rightarrow\infty}x_{n}=\ell$. Throughout this paper we assume that  $\liminf_{r}\frac{k_{r}}{k_{r-1}}>1$. $S_{\theta}(\textbf{R})$ will denote the set of lacunarily statistical convergent sequences of points in $\textbf{R}$. For an introduction to lacunary statistical convergence see \cite{FridyandOrhanlacunarystatisconvergence}.

Following the idea given in a 1946 American Mathematical Monthly problem by Buck \cite{Buck}, a number of authors Posner \cite{Posner}, Iwinski \cite{Iwinski}, Srinivasan \cite{Srinivasan}, Antoni \cite{Antoni}, Antoni and Salat \cite{AntoniSalat}, Spigel and Krupnik \cite{SpigelKrupnik} have studied $A$-continuity defined by a regular summability matrix $A$. Some authors, \"{O}zt\"{u}rk \cite{Ozturk}, Savas \cite{Savas}, Sava\c{s} and Das \cite{SavasDas}, Borsik, and Salat \cite{BorsikSalat} have studied $A$-continuity for methods of almost convergence or for related methods.

Recently, Connor and Grosse-Erdman \cite{ConnorandGrosse} have given sequential definitions of continuity calling $G$-continuity instead of $A$-continuity and their results cover the earlier works related to $A$-continuity where a method of sequential convergence, or briefly a method, is a linear function $G$ defined on a linear subspace of $s$, denoted by $c_{G}$, into $\textbf{R}$. A sequence $\textbf{x}=(x_{n})$ is said to be $G$-convergent to $\ell$ if $\textbf{x}\in c_{G}$ and $G(\textbf{x})=\ell$. In particular, $\lim$ denotes the limit function $\lim \textbf{x}=\lim_{n}x_{n}$ on the linear space $c$, $st-\lim$ denotes the statistical limit function $st-\lim \textbf{x}=st-\lim_{n}x_{n}$ on the linear space $S(\textbf{R})$ and $S_{\theta}-\lim$ denotes
the lacunary statistical limit function $S_{\theta}-\lim \textbf{x}=S_{\theta}-\lim_{n}x_{n}$ on the linear space $S_{\theta}(\textbf{R})$. A real function $f$ is called $G$-continuous at a point $u$ provided that whenever a sequence $\textbf{x}=(x_{n})$ of points in the domain of $f$ is $G$-convergent to $u$, then the sequence $f(\textbf{x})=(f(x_{n}))$ is $G$-convergent to $f(u)$. A method $G$ is called regular if every convergent sequence $\textbf{x}=(x_{n})$ is $G$-convergent with $G(\textbf{x})=\lim \textbf{x}$. A method is called subsequential if whenever $\textbf{x}$ is $G$-convergent with $G(\textbf{x})=\ell$,
then there is a subsequence $(x_{n_{k}})$ of $\textbf{x}$ with $\lim_{k} x_{n_{k}}=\ell$.

 A concept of quasi-Cauchy continuity, and a concept of quasi-Cauchy compactness (note that the terms forward continuity and forward compactness were used instead of quasi-Cauchy continuity and quasi-Cauchy compactness as well) have been introduced in \cite{CakalliForwardcontinuity} in the senses that a real function is called quasi-Cauchy continuous if $\lim_{n\rightarrow\infty}\Delta f(x_{n})=0$ whenever $\lim_{n\rightarrow\infty}\Delta x_{n}=0$, i.e. $f$ preserves quasi-Cauchy sequences; and a subset $E$ of $\textbf{R}$ is called quasi-Cauchy compact if any sequence $\textbf{x}=(x_{n})$ of points  in $E$   has a quasi-Cauchy subsequence (see also \cite{BurtonColeman}). Recently, some further results on quasi-Cauchy sequences are obtained in \cite{Cak3}, \cite{CakallistatisticalquasiCauchysequences}, and \cite{CakalliStatisticalwardcontinuity} (see also \cite{KizmazOncertainsequencespaces} and \cite{BasarirandSelmaAltundagDeltaLacunarystatistical}).
\maketitle

\section{$p$-quasi-Cauchy continuity}
We say that a sequence $\textbf{x}=(x_{n})$ is $p$-quasi-Cauchy with $\ell$ if $\lim_{k\rightarrow \infty} \Delta_{p} x_{k}=\ell$ where $\Delta_{p} x_{k}=x_{k+p}-x_{k}$. For the special case, $\ell=0$, $\textbf{x}$ is called $p$-quasi-Cauchy. We will deal with $p$ quasi-Cauchy sequences in this paper. Note that $\textbf{x}$ is quasi-Cauchy when $p=0$ and $\ell=0$, i.e. $1$-quasi-Cauchy sequences are quasi-Cauchy sequences. We will denote the set of all $p$-quasi-Cauchy sequences by $\Delta_{p}$. It follows from the equality $$x_{n+p}-x_{n}=x_{n+p}-x_{n+p-1}+x_{n+p-1}-x_{n+p-2}+x_{n+p-2}-+...+-x_{n+2}+x_{n+2}-x_{n+1}$$
 \;\;\;\;\;\;\;\;\;\;\;\;\;\;\;\;\;\;$+x_{n+1}-x_{n}$\\
that any quasi-Cauchy sequence is also $p$-quasi-Cauchy, but the converse is not always true as it can be seen by considering the sequence $((-1)^{n})$. Any any slowly oscillating sequence is $p$-quasi-Cauchy, so is Cauchy sequence. But neither, being $\delta$-quasi-Cauchy implies $p$-quasi-Cauchy, nor being being $p$-quasi-Cauchy implies being $\delta$-quasi-Cauchy. Counterexamples for these situations are the sequences $(n)$, and $((-1)^{n})$. We note that the sum of two $p$-quasi-Cauchy sequence is $p$-quasi-Cauchy whereas product of two $p$-quasi-Cauchy sequences need not be $p$-quasi-Cauchy as it can be seen by considering the product of the sequence $(\sqrt{n})$ itself.

Now we give the definition of $p$-quasi-Cauchy compactness of a subset of $\textbf{R}$.

\begin{Def}
A subset $E$ of $\textbf{R}$ is called $p$-quasi-Cauchy compact if whenever $\textbf{x}=(x_{n})$ is a sequence of points in $E$ there is a $p$-quasi-Cauchy subsequence $\textbf{z}=(y_{k})=(x_{n_{k}})$ of $\textbf{x}$.
\end{Def}

We note that this definition of $p$-quasi-Cauchy compactness cannot be obtained by any summability matrix $A$, even by the summability matrix $A=(a_{kn})$ defined by $a_{kn}=-1$ if $k=n$, and $a_{kn}=1$ if $k=n+p$; $a_{kn}=0$ otherwise,

\[G(x)=\lim A\textbf{x}=\lim_{k\rightarrow\infty}\sum^{\infty}_{n=1}a_{kn}a_{n}=\lim_{k\rightarrow \infty} \Delta_{p} x_{k}  \;  \;  \;  \; \; (*)
\]
(see \cite{CakalliSequentialdefinitionsofcompactness} for the definition of $G$-sequential compactness).
In this example, despite that $G$-sequential compact subsets of $\textbf{R}$ should include the singleton set $\{0\}$, $p$-quasi-Cauchy compact subsets of $\textbf{R}$ do not have to include the singleton $\{0\}$.

Since any quasi-Cauchy sequence is $p$-quasi-Cauchy we see that any ward compact (quasi-Cauchy compact) subset of $\textbf{R}$ is $p$-quasi-Cauchy compact for any $p\in{\textbf{N}}$. We see that any finite subset of $\textbf{R}$ is $p$-quasi-Cauchy compact, union of two $p$-quasi-Cauchy compact subsets of $\textbf{R}$ is $p$-quasi-Cauchy compact and intersection of any $p$-quasi-Cauchy compact subsets of $\textbf{R}$ is $p$-quasi-Cauchy compact. Furthermore any subset of a $p$-quasi-Cauchy compact set of \textbf{R} is $p$-quasi-Cauchy compact and any bounded subset of $\textbf{R}$ is $p$-quasi-Cauchy compact. Any compact subset of $\textbf{R}$ is also $p$-quasi-Cauchy compact and the converse is not always true. For example, the set $\{1, \frac{1}{2}, \frac{1}{3},...,\frac{1}{n},...\}$ is $p$-quasi-Cauchy compact for every $p\in{\textbf{N}}$, which is not compact. On the other hand, the set $\textbf{N}$ is not $p$-quasi-Cauchy compact for any $p\in{\textbf{N}}$. We note that any slowly oscillating compact subset of $\textbf{R}$ is $p$-quasi-Cauchy compact (see \cite{CakalliSlowlyoscillatingcontinuity}  for the results on slowly oscillating compactness). These observations above suggest to us the following:

\begin{Lem}
 A subset $A$ of \textbf{R} is bounded if and only if it is p-quasi-Cauchy compact.
\end{Lem}
\begin{pf}
It is to see that bounded subsets of \textbf{R} are $p$-quasi-Cauchy compact. To prove the converse suppose that $A$ is unbounded. If it is unbounded above, then one can construct a sequence $(x _{n})$ of numbers in  $A$ such that  $x_{n+1}>p+x _{n}$ for each positive integer $n$. Here we note that $x_{n+1}>1+x _{n}$ will suffice. Then the sequence $(x_{n})$ does not have any $p$-quasi-Cauchy subsequence, so $A$ is not $p$-quasi-Cauchy compact. If $A$ is bounded above and unbounded below, then similarly we obtain that $A$ is not $p$-quasi-Cauchy compact. This completes the proof.
\end{pf}

\begin{Cor}
\normalfont
A subset $A$ of \textbf{R} is bounded if and only if any sequence of points in $A$ has a subsequence which is any type of the following: Cauchy, quasi-Cauchy, $p$-quasi-Cauchy for any $p\in{\textbf{N}}$, statistically quasi-Cauchy, lacunary statistically quasi-Cauchy, and slowly oscillating.
\end{Cor}
\begin{pf}
The proof follows from the preceding theorem and Theorem 3 in \cite{CakallistatisticalquasiCauchysequences} easily so is omitted.
\end{pf}

The concept of sequential continuity suggests to us giving a new type continuity, namely, $p$-quasi-Cauchy continuity, analogous to quasi-Cauchy continuity.

\begin{Def}
A function $f$ is called $p$-quasi-Cauchy continuous on $E$ if the sequence $f(\textbf{x})=(f(x_{n}))$ is $p$-quasi-Cauchy whenever $\textbf{x}=(x_{n})$ is a sequence of terms in $E$ which is $p$-quasi-Cauchy, i.e. $f$ preserves $p$-quasi-Cauchy sequences.
\end{Def}
We note that this definition of $p$-quasi-Cauchy continuity can not be obtained by any summability matrix $A$, even by the summability matrix $A=(a_{nk})$ defined by $(*)$ however for this special summability matrix $A$ if $A$-continuity of $f$ at the point $0$ implies $p$-quasi-Cauchy continuity of $f$, then $f(0)=0$; and if $p$-quasi-Cauchy  continuity of $f$ implies $A$-continuity of $f$ at the point $0$, then $f(0)=0$.

We also note that sum of two $p$-quasi-Cauchy continuous functions is
$p$-quasi-Cauchy continuous, but the product of
two $p$-quasi-Cauchy continuous functions need not be $p$-quasi-Cauchy continuous as it can be seen by considering product of the
$p$-quasi-Cauchy continuous function $f(x)=x$ with itself.

In connection with $p$-quasi-Cauchy sequences and convergent sequences the problem arises to investigate the following types of  "continuity" of a function on $\textbf{R}$.

\begin{description}
\item[($\Delta_{p} $)] $(x_{n}) \in {\Delta_{p}} \Rightarrow (f(x_{n})) \in {\Delta_{p}}$
\item[($\Delta_{p} c$)] $(x_{n}) \in {\Delta_{p}} \Rightarrow (f(x_{n})) \in {c}$
\item[$(c)$] $(x_{n}) \in {c} \Rightarrow (f(x_{n})) \in {c}$
\item[$(d)$] $(x_{n}) \in {c} \Rightarrow (f(x_{n})) \in {\Delta_{p}}$
\end{description}

We see that $(\Delta_{p} )$ is $p$-quasi-Cauchy continuity of $f$, and $(c)$ states the ordinary continuity of $f$. It is easy to see that $(\Delta_{p} c)$ implies $(\Delta_{p})$, and $(\Delta_{p} )$ does not imply $(\Delta_{p} c)$;  and $(\Delta_{p})$ implies $(d)$, and $(d)$ does not imply $(\Delta_{p}c)$; $(\Delta_{p} c)$ implies $(c)$ and $(c)$ does not imply $(\Delta_{p} c)$; and $(c)$ is equivalent to $(d)$.

There is no function which satisfies that $  (x_{n}) \in {c} \Rightarrow (f(x_{n}))\in {c_{0}}$ where $c_{0}$ denotes the set of all null sequences.
We note that if $  (x_{n}) \in {c_{0}}$ implies that $(f(x_{n}))\in {c}$, then $\lim_{n\rightarrow\infty} f(x_{n})=f(0)$.

Now we give the implication $(\Delta_{p})$ implies $(\Delta_{1})$, i.e. any $p$-quasi-Cauchy continuous function is $1$-quasi-Cauchy continuous, i.e. quasi-Cauchy continuous.

\begin{Thm} \label{TheopquasiCauchyimpliesquasiCauchy} If $f$ is $p$-quasi-Cauchy continuous on a subset $E$ of $\textbf{R}$ for a $p\in{\textbf{N}}$, then it is quasi-Cauchy continuous
on $E$.
\end{Thm}
\begin{pf} Let $p$ be any positive integer. If $p=1$, then there is nothing to prove. So we would suppose that $p>1$. Take any $p$-quasi-Cauchy continuous function $f$ on $E$. Let $(x_{n})$ be any quasi-Cauchy sequence of points in $E$. Then the sequence $$(x_{1}, x_{1},..., x_{1}, x_{2}, x_{2},..., x_{2},..., x_{n}, x_{n},..., x_{n}, ...)$$ is also quasi-Cauchy so it is $p$-quasi-Cauchy, hence
it belongs to $\Delta_{p}$  so does the sequence
$$(f(x_{1}), f(x_{1}), ... ,f(x_{1}),f(x_{2}), f(x_{2}),..., f(x_{2}),..., f(x_{n}), f(x_{n}), ..., f(x_{n}), ...)$$ where the same term repeats $p$-times. Thus $\lim_{n\rightarrow\infty}|f(x_{n+1})-f(x_{n})|=0$.
This completes the proof of the theorem.

\end{pf}

\begin{Cor} If $f$ is $p$-quasi-Cauchy continuous, then it is continuous in the ordinary case.
\end{Cor}
\begin{pf} The proof follows immediately from Theorem 1 on page 228 in \cite{CakalliForwardcontinuity} so is omitted.
\end{pf}

The concept of slowly oscillating continuity was introduced by Cakalli in \cite{CakalliSlowlyoscillatingcontinuity}, and investigated by Canak and Dik \cite{DikandCanak}. Vallin also studied slowly oscillating continuity in \cite{Vallin}, but his approach to the proof that any slowly oscillating continuous function is continuous is an indirect way although Cakalli's proof is direct, easy, and nice depending only on the fact that any slowly oscillating sequence is quasi-Cauchy. Thus we have the following.

\begin{Thm}  $p$-quasi-Cauchy continuous image of any $p$-quasi-Cauchy compact subset of $\textbf{R}$ is $p$-quasi-Cauchy compact.
\end{Thm}
\begin{pf}
Let $f$ be a $p$-quasi-Cauchy continuous function and $E$ be a $p$-quasi-Cauchy compact subset of $\textbf{R}$. Take any sequence $\textbf{y}=(y_{n})$ of terms in $f(E)$. Write $y_{n}=f(x_{n})$ where $x_{n}\in {E}$ for each $n \in{\textbf{N}}$. $p$-quasi-Cauchy compactness of $E$ implies that there is a subsequence
$\textbf{y}=(y_{k})=(x_{n_{k}})$ of $\textbf{x}$ with $\lim_{k\rightarrow \infty} \Delta_{p}y_{k}=0$. Since $f$ is $p$-quasi-Cauchy continuous, $(t_{k})=f(\textbf{y})=(f(y_{k}))$ is $p$-quasi-Cauchy . Thus $(t_{k})$ is a subsequence of the sequence $f(\textbf{x})$ with $\lim_{k\rightarrow \infty} \Delta_{p} t_{k}=0$. This completes the proof of the theorem.
\end{pf}

\begin{Cor} $p$-quasi-Cauchy continuous image of any compact subset of $\textbf{R}$ is compact.
\end{Cor}
\begin{pf} The proof follows from the preceding theorem, so is omitted.
\end{pf}

Now we prove that any uniformly continuous function preserves $p$-quasi-Cauchy sequences.

\begin{Thm}  \label{TheouniformcontibutyimpliþespQuasiCauchycontinuity} If $f$ is uniformly continuous on a subset $E$ of $\textbf{R}$, then it is $p$-quasi-Cauchy continuous on $E$ for any $p\in{\textbf{N}}$.
\end{Thm}
\begin{pf}
To prove that $(f(x_{n}))$ is a $p$-quasi-Cauchy sequence whenever $(x_{n})$ is, take any $\varepsilon > 0$. Uniform continuity of $f$ on $E$ implies that there exists a $\delta >0$, depending on $\varepsilon$, such that $|f(x)-f(y)|< \varepsilon$ whenever  $|x-y|< \delta$ and $x, y\in{E}$. For this $\delta >0$, there exists an $N=N(\delta)=N_{1}(\varepsilon)$ such that $|\Delta_{p} x_{n}|<\delta$ whenever $n>N$. Hence $|\Delta_{p} f(x_{n})|<\varepsilon$ if $n>N$. It follows from this that $(f(x_{n}))$ is a $p$-quasi-Cauchy sequence. This completes the proof of the theorem.
\end{pf}

\begin{Cor} If $f$ is slowly oscillating continuous on an interval $E$, then it is $p$-quasi-Cauchy continuous for all $p\in{\textbf{N}}$.
\end{Cor}
\begin{pf}
If $f$ is a slowly oscillating continuous function on an interval $E$, then it is uniformly continuous on $E$ by Theorem 5 on page 1623 in \cite{CakallistatisticalquasiCauchysequences}. Hence it follows from the preceding theorem that $f$ is $p$-quasi-Cauchy continuous on $E$ for all $p\in{\textbf{N}}$.
\end{pf}

It is well-known that any continuous function on a compact subset $E$ of $\textbf{R}$ is uniformly continuous on $E$.
We have an analogous theorem for a $p$-quasi-Cauchy continuous function defined on a $p$-quasi-Cauchy compact subset of $\textbf{R}$.

\begin{Thm} If a function is $p$-quasi-Cauchy continuous on a $p$-quasi-Cauchy compact subset $E$ of $\textbf{R}$, then it is uniformly continuous on $E$.
\end{Thm}

\begin{pf} Suppose that $f$ is not uniformly continuous on $E$ so that there exist an $\epsilon _{0}>0$ and sequences $(x_{n})$ and $(y_{n})$ of points in $E$ such that
\[
|x_{n}-y_{n}|<1/n
\]
and
\[
|f(x_{n})-f(y_{n})|\geq \epsilon _{0}
\]
for all $n \in \textbf{N}$.
Since $E$ is $p$-quasi-Cauchy compact, there is a subsequence $(a_{n_{k}})$ of $(x_{n})$ that is $p$-quasi-Cauchy. On the other hand there is a subsequence $(y_{n_{k_{j}}})$ of $(y_{n_{k}})$ that is $p$-quasi-Cauchy as well. It is clear that the corresponding sequence $(a_{n_{k_{j}}})$ is also $p$-quasi-Cauchy, since $(y_{n_{k_{j}}})$ is $p$-quasi-Cauchy and
\[|x_{n_{k_{j}}}-x_{n_{k_{j+p}}}|\leq |x_{n_{k_{j}}}-x_{n_{k_{j}}}|+|x_{n_{k_{j}}}-x_{n_{k_{j+p}}}|+|x_{n_{k_{j+p}}}-x_{n_{k_{j+p}}}|.
\]
which follows from the inequality
\[|y_{n_{k_{j+p}}}-x_{n_{k_{j+p}}}|\leq |y_{n_{k_{j+p}}}-y_{n_{k_{j}}}|+|y_{n_{k_{j}}}-x_{n_{k_{j}}}|+|x_{n_{k_{j}}}-x_{n_{k_{j+p}}}|.
\]
Now the sequence $$(x _{n_{k_{1}}}, x_{n_{k_{1}}},..., x_{n_{k_{1}}}, y_{n_{k_{1}}}, y_{n_{k_{1}}},... ,y_{n_{k_{1}}},  ..., x_{n_{k_{j}}}, x_{n_{k_{j}}},...,x_{n_{k_{j}}}, y_{n_{k_{j}}}, y_{n_{k_{j}}}, ..., y_{n_{k_{j}}},...)$$
is $p$-quasi-Cauchy while the transformed sequence is not $p$-quasi-Cauchy where same terms repeat $p$-times. Hence this establishes a contradiction so completes the proof of the theorem.
\end{pf}

We see that for a regular subsequential method $G$ that any $p$-quasi-Cauchy continuous function on a $G$-sequentially compact subset $E$ of $\textbf{R}$ is uniformly continuous on $E$ (see \cite{CakalliSequentialdefinitionsofcompactness} and \cite{CakalliOnGcontinuity}).

\begin{Lem} \label{Lemorderedpairsquasi-Cauchy} If $(\xi_{n}, \eta_{n})$ is a sequence of ordered pairs of points in an interval $E$ of \textbf{R} such that $\lim_{n\rightarrow\infty} |\xi_{n}-\eta_{n})|=0$, then there exists a $p$-quasi-Cauchy sequence $(x_{n})$ with the property that for any positive integer $i$ there exists a positive integer $j$ such that $(\xi_{i}, \eta_{i})=(x_{j-p}, x_{j})$.
\end{Lem}
\begin{pf}
For each positive integer $k$, we can fix $y_{0}^{k}, y_{1}^{k}, ...,y_{n_{k}}^{k}$ in $E$ with $y_{0}^{k}=\eta_{k}$, $y_{n_{k}}^{k}=\xi_{k+1}$, and $|y_{i}^{k}- y_{i-p}^{k}|<\frac{1}{k}$  for $1\leq{i}\leq{n_{k}}$. Now write

$$ (\xi_{1}, \eta_{1}, y_{1}^{1},...,y_{n_{1}-1}^{1}, \xi_{2}, \eta_{2}, y_{1}^{2},...,y_{n_{2}-1}^{2}, \xi_{3}, \eta_{3},..., \xi_{k}, \eta_{k}, y_{1}^{k}, ...,y_{n_{k-1}}^{k}, \xi_{k+1}, \eta_{k+1},...)$$
Then denoting this sequence by $(x_{n})$ we obtain that for any positive integer $i$ there exists a positive integer $j$ such that $(\xi_{i}, \eta_{i})=(x_{j-p}, x_{j})$. This completes the proof of the lemma.

\end{pf}

\begin{Thm}  If $f$ is is $p$-quasi-Cauchy continuous on an interval $E$ for a positive integer $p$, then it is uniformly continuous on $E$.
\end{Thm}
\begin{pf}
Let $p$ be any positive integer and $E$ an interval. To prove that $p$-quasi-Cauchy continuity of $f$ on $E$ implies uniform continuity of $f$ on $E$ suppose that $f$ is not uniformly continuous on $E$ so that there exists an  $\varepsilon > 0$ such that for any $\delta >0$, there exist $x, y \in{E}$ with $|x-y|<\delta$ but $|f(x)-f(y)|\geq \varepsilon$. Hence for each positive integer $n$, there exist $x_n$ and $y_n$ in $E$ such that $|x_{n}-y_{n})<\frac{1}{n}$, and $|f(x_{n})-f(y_{n})|\geq \varepsilon$. By Lemma \ref{Lemorderedpairsquasi-Cauchy}, one can construct a $p$-quasi-Cauchy sequence  $(t_{n})$ which has a subsequence $(z_{n})=(t_{k_{n}})$ such that  $\lim_{n\rightarrow\infty} ( z_{n+p}-z_{n})=0$, but  $|f(z_{n+p})-f(z_{n})|\geq \varepsilon$. Therefore the transformed sequence $(f(z_{n}))$ is not $p$-quasi-Cauchy  Thus this contradiction yields that $p$-quasi-Cauchy continuity implies uniform continuity. This completes the proof of the theorem.

\end{pf}

Combining Theorem \ref{TheouniformcontibutyimpliþespQuasiCauchycontinuity} and the preceding theorem we have the following:

\begin{Cor} A function $f$ is uniformly continuous on an interval $E$ if and only if there exists a positive integer $p$ such that $f$ preserves $p$-quasi-Cauchy sequences of points in $E$.

\end{Cor}

\begin{Cor} A function defined on an interval is $p$-quasi-Cauchy continuous for a $p\in{\textbf{N}}$ if and only if it is slowly oscillating continuous.
\end{Cor}
\begin{Cor} If $f$ is $\delta$-quasi-Cauchy continuous on a subset $E$ of $\textbf{R}$, then it is $p$-quasi-Cauchy continuous on $E$.
\end{Cor}
\begin{pf}
It follows from Theorem 7 on page 399 in \cite{CakalliDeltaquasiCauchysequences}.
\end{pf}

It is a well known result that uniform limit of a sequence of continuous functions is continuous. This is also true in the case of $p$-quasi-Cauchy continuity; i.e. uniform limit of a sequence of $p$-quasi-Cauchy continuous functions is $p$-quasi-Cauchy continuous.

\begin{Thm} If $(f_{n})$ is a sequence of $p$-quasi-Cauchy continuous functions defined on a subset $E$ of $\textbf{R}$ and $(f_{n})$ is uniformly convergent to a function $f$, then $f$ is $p$-quasi-Cauchy continuous on $E$.
\end{Thm}
\begin{pf} Let $\textbf{x}=(x_{n})$ be any $p$-quasi-Cauchy  sequence of points in $E$ and $\varepsilon > 0$. Then there exists a positive integer $N$ such that $|f_{n}(x)-f(x)|<\frac{\varepsilon}{3}$ for all $x \in {E}$ whenever $n\geq N$. As $f_{N}$ is $p$-quasi-Cauchy continuous, there exists a positive integer $N_{1}$, depending on $\varepsilon$ and greater than $N$ such that $|f_{N}(x_{n+p})-f_{N}(x_{n})|<\frac{\varepsilon}{3}$ for $n\geq N_{1}$. Now for $n\geq N_{1}$  we have \\$ |f(x_{n+p})-f(x_{n})|
 \leq |f(x_{n+p})-f_{N}(x_{n+p})|+|f_{N}(x_{n+p})-f_{N}(x_{n})|+|f_{N}(x_{n})-f(x_{n})|$
\; \; \; \; \;  $< \frac{\varepsilon}{3} + \frac{\varepsilon}{3} + \frac{\varepsilon}{3}=\varepsilon.$\\
This completes the proof of the theorem.
\end{pf}

\begin{Thm}  The set of all $p$-quasi-Cauchy continuous functions defined on a subset $E$ of $\textbf{R}$ is a closed subset of the set of all continuous
functions on $E$, i.e. $\overline{\Delta_{p}FC(E)}=\Delta_{p}FC(E)$ where $\Delta_{p}FC(E)$ is the
set of all $p$-quasi-Cauchy continuous functions on $E$, $\overline{\Delta_{p}FC(E)}$ denotes
the set of all cluster points of $\Delta_{p}FC(E)$ and $E$ is a bounded subset
of $\textbf{R}$.
\end{Thm}

\begin{pf} Let us denote the set of all $p$-quasi-Cauchy continuous functions on $E$ by $\Delta_{p}FC(E)$ and $f$ be any element in the closure of $\Delta_{p}FC(E)$. Then there exists a sequence of points in $\Delta_{p}FC(E)$ such that $\lim_{k\rightarrow \infty} f_{k}=f$.
To show that $f$ is $p$-quasi-Cauchy continuous take any $p$-quasi-Cauchy  sequence $(x_{n})$. Let $\varepsilon > 0$. Since $(f_{k})$
converges to $f$, there exists an $N$ such that for all $x \in {E}$
and for all $n \geq {N}$, $|f(x)-f_{n}(x)|< \frac{\varepsilon}{3}$.
As $f_{N}$ is $p$-quasi-Cauchy continuous, there is an $N_{1}$, greater
than $N$, such that for all $n \geq {N_{1}}$,
$|f_{N}(x_{n+p})-f_{N}(x_{n})|<
\frac{\varepsilon}{3}$. Hence for all $n \geq {N_{1}}$,
\\$|f(x_{n+p})-f(x_{n})|$ $\leq {|f(x_{n+p})-f_{N}(x_{n+p})|+|f(x_{n})-f_{N}(x_{n})|}$ $+|f_{N}(x_{n+p})-f_{N}(x_{n})|$ \; \; \; \; \;
 $< \frac{\varepsilon}{3} + \frac{\varepsilon}{3} + \frac{\varepsilon}{3}= \varepsilon.$\\
This completes the proof of the theorem.
\end{pf}

\begin{Cor} The set of all $p$-quasi-Cauchy continuous functions on a subset $E$ of $\textbf{R}$ is a complete subspace of
the space of all continuous functions on $E$.
\end{Cor}
\begin{pf}
The proof follows from the preceding theorem.
\end{pf}

\section{Conclusion}
In this paper we introduce a concept of $p$-quasi-Cauchy continuity via $p$-quasi-Cauchy sequences and investigate results related to this kind of continuity and the some other kinds of continuities; namely ordinary continuity, uniform continuity, statistical continuity, slowly oscillating continuity, and $\delta$-quasi-Cauchy continuity. For the special case $p=1$ we obtain most of the results in \cite{CakalliForwardcontinuity}, \cite{BurtonColeman}, \cite{Cak3}, and \cite{CakalliStatisticalwardcontinuity}.

For a further study, we suggest to investigate $p$-quasi-Cauchy sequences of fuzzy points and $p$-quasi-Cauchy continuity for the fuzzy functions. However due to the change in settings, the definitions and methods of proofs will not always be analogous to those of the present work (for example see \cite{CakalliandPratul}). For another further study we also suggest to introduce a new concept in dynamical systems using $p$-quasi-Cauchy continuity.

\end{document}